\documentclass[12pt]{article}%
\usepackage{amsmath}
\usepackage{amsfonts}
\usepackage{amssymb}
\usepackage{graphicx}%
\providecommand{\U}[1]{\protect\rule{.1in}{.1in}}
\newtheorem{theorem}{Theorem}

\newtheorem{corollary}[theorem]{Corollary}

\newtheorem{example}[theorem]{Example}

\newtheorem{lemma}[theorem]{Lemma}

\newtheorem{proposition}[theorem]{Proposition}
\newtheorem{remark}[theorem]{Remark}

\newenvironment{proof}[1][Proof]{\noindent\textbf{#1.} }{\ \rule{0.5em}{0.5em}}
\begin{document}

\title{Efficient orders for reciprocal matrices with block structure}
\author{Susana Furtado\thanks{Corresponding author. Email: sbf@fep.up.pt.
orcid.org/0000-0003-0395-5972. The work of this author is funded by FCT -
Funda\c{c}\~{a}o para a Ci\^{e}ncia e a Tecnologia, I.P., through national
funds, under the project UID/04561/2025,
https://doi.org/10.54499/UID/04561/2025.}\\CEMS.UL and Faculdade de Economia e Gest\~{a}o\\Universidade do Porto\\Rua Dr. Roberto Frias\\4200-464 Porto, Portugal
\and Charles R. Johnson \thanks{Email: crjmatrix@gmail.com. }\\225 West Tazewells Way\\Williamsburg, VA 23185, United States}
\maketitle

\begin{abstract}
We consider block partitioned reciprocal matrices. For these matrices, we
characterize the efficient vectors whose subvectors, corresponding to
principal blocks, are positive multiples of given efficient vectors for those
blocks. This generalizes a prior, scalar oriented, extension theorem.
Conditions are presented that guarantee the existence of efficient vectors
that preserve ordinal relations associated with cross-block comparisons. The
notion of order transitivity is also generalized to block reciprocal matrices.
For such matrices, it is then shown that efficient vectors, with no rank
reversals in positions corresponding to off-diagonal blocks, exist. The notion
of uniformly ordered reciprocal matrices (all efficent vectors exhibit the
same order) is also generalized in the block reciprocal setting. This class is
characterized via known matrix bounds for the consistent matrices associated
with efficient vectors. We also provide a graph-theoretic characterization of
the existence of a positive vector that is efficient for all completions of a
given partial reciprocal matrix. Based on this, we apply some of our results
on ordering to such matrices.

\end{abstract}


\bigskip

\section{Introduction}

An $n$-by-$n$ entrywise positive matrix $A=[a_{ij}]$ is called
\emph{reciprocal} if $a_{ji}=\frac{1}{a_{ij}}$ for all $1\leq i,j\leq n$. The
diagonal entries are all $1$, and the off-diagonal entries represent pairwise
ratio comparisons among $n$ alternatives. We denote by $\mathcal{PC}_{n}$ the
set of all such matrices. If $A\in\mathcal{PC}_{n}$ satisfies, in addition,
$a_{ij}a_{jk}=a_{ik}$, $1\leq i,j,k\leq n$, then $A$ is said to be
\emph{consistent}. Any matrix in $\mathcal{PC}_{n}$, $n\leq2,$ is consistent.

In a variety of ranking models, as the Analytic Hierarchy Process
\cite{saaty1977}, a weight vector $w=\left[
\begin{array}
[c]{ccc}%
w_{1} & \cdots & w_{n}%
\end{array}
\right]  ^{T}\in\mathbb{R}_{+}^{n}$ (the set of positive $n$-vectors) is to be
deduced from $A\in\mathcal{PC}_{n}$, so that $\frac{w_{i}}{w_{j}}\approx
a_{ij}$. If $A$ is consistent, then $A$ is generated by a (projectively)
single positive vector $w\in\mathbb{R}_{+}^{n}$, so that $A=ww^{(-T)}=\left[
\frac{w_{i}}{w_{j}}\right]  $. Here, $w^{(-T)}$ is the entrywise inverse of
the transpose of $w$. In this case, $w$ is (projectively) the natural weight
vector. However, consistency of comparisons seldom occurs in practice. For
this reason, it is important to obtain weight vectors for an inconsistent
reciprocal matrix. We denote the set of inconsistent matrices in
$\mathcal{PC}_{n}$ by $\mathcal{PC}_{n}^{0}$. Several methods for obtaining
weight vectors have been proposed and discussed, see surveys
\cite{Choo,Golany1993}.

We ask that the consistent matrix $W=ww^{(-T)}$ approximates $A\in
\mathcal{PC}_{n}$ in such a way that no other $v\in\mathbb{R}_{+}^{n}$ better
approximates $A$ via $V=vv^{(-T)}$, i.e. $\left\vert A-V\right\vert
\leq\left\vert A-W\right\vert ,$ entrywise, implies that $v$ is proportional
to $w.$ Such a vector $w$ is called \emph{efficient} (or \emph{Pareto
optimal}). We denote the set of all efficient vectors for $A$ by
$\mathcal{E}(A)$, which will be projectively infinite if $A\in\mathcal{PC}%
_{n}^{0}$. For a consistent matrix, $\mathcal{E}(A)$ is projectively a single
vector. Two methods to describe all efficient vectors for each $A\in
\mathcal{PC}_{n}^{0}$ have recently been developed; $\mathcal{E}(A)$ as a
union of a finite number of convex sets \cite{FJ3}, and an inductive
construction of $\mathcal{E}(A)$ \cite{FJ4}. The former will be used here and
is described in Section \ref{sunion}. The latter gives a method to construct
$\mathcal{E}(A)$ from the efficient vectors for the principal submatrices. It
also happens that any efficient vector for a principal submatrix of $A$ can be
extended to a vector in $\mathcal{E}(A)$ \cite{FJ1}. This will be generalized.

In many applications, the goal is to order the alternatives based on the
chosen weight vector. In these situations, more important than the cardinal
ranking is an ordinal one. This motivates the study of the possible orders
that occur among efficient vectors, called \emph{efficient orders}. Unless
otherwise stated, by an order of a vector we mean a weak order. If
$w\in\mathbb{R}_{+}^{n}$ with $w_{i}\geq w_{i+1},$ for $i=1,\ldots,n-1$, we
say that $w$ has decreasing order.

We say that $w\in\mathbb{R}_{+}^{n}$ has a \emph{rank reversal} for
$A\in\mathcal{PC}_{n}$ at position $i,j$ if one of the following holds:
$a_{ij}>1$ and $\frac{w_{i}}{w_{j}}<1$, or $a_{ij}<1$ and $\frac{w_{i}}{w_{j}%
}>1.$ This means that, for alternatives $i$ and $j$, the weight vector $w$ is
order inconsistent with the pairwise comparison $a_{ij}$. In $A$, alternative
$i$ is judged more (less) important than alternative $j$ but, in $w$, $i$ is
assigned a lower (higher) weight than $j$. Note that, if $w$ has a rank
reversal for $A$ at position $i,j$, then it also has one at position $j,i$. In
general, an efficient order minimizing the number of rank reversals is
desirable. In \cite{FJRR} this question was studied in detail.

An efficient vector preserving the ordinal judgments given by the reciprocal
matrix along a labelled spanning tree always exists. In fact, a stronger
statement holds: there is an efficient vector for $A$ such that the
corresponding consistent matrix coincides with $A$ along the spanning tree.
The reciprocal matrices for which there is at least one weight vector that has
no rank reversals for $A$ are called \emph{order transitive}. For such
matrices, the existence of an efficient weight vector with no rank reversals
was shown in \cite{FJRR}. In \cite{FJ20}, the matrices in $\mathcal{PC}_{n}$
for which all efficient vectors have the same order were studied for the first
time. These were called \emph{uniformly ordered matrices} and are precisely
the reciprocal matrices for which no efficient vector gives a rank reversal.
Several characterizations of these matrices were presented in \cite{FJ8,FJRR}.

Here we generalize the above mentioned results to reciprocal matrices
partitioned into blocks. (The partition of the row indices is the same as the
one of the column indices.) We refer to such matrices as \emph{block
reciprocal matrices}.

\bigskip The paper is organized as follows. We start by introducing some
necessary background in Section \ref{sbackground}. Then, in Section
\ref{saddobs}, besides some simple observations, we give necessary and
sufficient conditions for every ordering to be realized in a constituent
convex set of $\mathcal{E}(A).$ For incomplete reciprocal matrices, we also
provide a graph-theoretic characterization of the existence of a positive
vector that is efficient for all completions. In Section \ref{Seffprincsubm},
we consider block reciprocal matrices $A$ and characterize the vectors in
$\mathcal{E}(A)$ whose subvector corresponding to each principal block is a
positive multiple of a prescribed efficient vector for that block (hence the
same consistent matrix is determined). There always exists such a vector. This
result generalizes the extension result in \cite{FJ1}. Section
\ref{Spreservord} is devoted to the preservation of ordinal information. We
prove that, given a partition of $\{1,\ldots,n\}$ into $p$ subsets and a
labelled spanning tree of the complete graph on $p$ vertices, if for any two
subsets, corresponding to an edge in this tree, all cross-sets comparisons
induce the same ordinal judgment, then there exists an efficient vector
preserving the corresponding ordering. Repeated applications of this result
allow the analysis of more general comparison patterns. We also apply the
result to partial reciprocal matrices. In Section \ref{SOT} we extend the
concept of order transitivity to block reciprocal matrices. We establish a
structural characterization of block reciprocal transitive matrices, showing
that they are block permutationally similar to block matrices whose
off-diagonal blocks above the diagonal contain no entries smaller than one. We
further prove the existence of efficient vectors exhibiting no rank reversals
at the corresponding block positions, thereby extending known results for
standard order transitive matrices. Finally, in Section \ref{SUO}, we
generalize the notion of uniformly ordered matrices to the setting of block
reciprocal matrices. We provide a characterization of this class in terms of
known matrix bounds for consistent matrices generated by efficient vectors.
Section \ref{sconcl} summarizes the main conclusions. Throughout the paper,
examples illustrating the theoretical results are presented.

\section{Background\label{sbackground}}

\subsection{Notation and basic observations}

For an integer $n>0,$ we denote by $N$ the set $\{1,\ldots,n\}$. If $A$ is an
$n$-by-$n$ matrix, $w\in\mathbb{R}_{+}^{n}$ and $S_{1},S_{2}\subset N,$ we
denote by $A[S_{1},S_{2}]$ the submatrix of $A$ with rows indexed by $S_{1}$
and columns indexed by $S_{2}.$ We simply write $A[S_{1},S_{1}]$ as
$A[S_{1}].$ By $w[S_{1}]$ we denote the subvector of $w$ with entries indexed
by $S_{1}$ and by $w_{i}$ the $i$th entry of $w$. For a positive matrix $Q$,
we denote by $Q^{(-T)}$ the transpose of its entrywise inverse. We denote by
$\mathbf{e}$ the vector of ones of appropriate size, and by $J$ a matrix all
of whose columns are $\mathbf{e}$. Given $w\in\mathbb{R}_{+}^{n}$, we often
denote the consistent matrix $ww^{(-T)}$ by $W$.

\bigskip

A cycle $\alpha$ in a directed graph (digraph) is a path (no repeated
vertices) with an additional edge from the last to the first vertex:
$i_{1}\rightarrow i_{2}\rightarrow i_{3}\rightarrow\cdots\rightarrow
i_{k}\rightarrow i_{k+1}=i_{1}$ (or $i_{1}i_{2}i_{3}\cdots i_{k}i_{1}$).
A\ path is called \emph{Hamiltonian} (an H-path) if it contains all vertices.
If a cycle contains an H-path, it is called a \emph{Hamiltonian cycle} (an H-cycle).

A directed graph is said to be \emph{semicomplete} if there is an edge in at
least one direction between any two distinct vertices. The following result
extends R\'{e}dei's Theorem for tournaments.

\begin{theorem}
\cite{Jensen,Moon}\label{thRedeis} Every semicomplete digraph contains a
Hamiltonian path.
\end{theorem}

Let $A\in\mathcal{PC}_{n}$ and $H$ be the digraph with vertex set $N$ and an
edge $(i,j)$ if and only if $a_{ij}\geq1.$ Then $H$ is semicomplete. By
Theorem \ref{thRedeis}, $H$ contains an H-path. This implies that there is an
$A^{\prime}\in\mathcal{PC}_{n}$, permutationally similar to $A$, such that the
digraph $H^{\prime}$ associated with $A^{\prime}$ contains the H-path
$12\cdots n.$ Thus, the following holds. By superdiagonal entries we mean the
entries in positions $i,i+1,$ for $i=1,\ldots,n-1.$

\begin{theorem}
\label{Prop1}Any reciprocal matrix is permutationally similar to a matrix with
all the superdiagonal entries $\geq1$.
\end{theorem}

This result was obtained in \cite{FJRR}, using different arguments, when all
the off-diagonal entries of the reciprocal matrix are distinct from $1.$ The
given proof also proves R\'{e}dei's Theorem for tournaments.

\subsection{Graph characterization of efficiency}

A \emph{monomial matrix} is a product of a permutation matrix and a positive
diagonal matrix.

\begin{lemma}
\label{lsim}\cite{FJ1} If $A\in\mathcal{PC}_{n}$ and $S$ is an $n$-by-$n$
monomial matrix, then $SAS^{-1}\in\mathcal{PC}_{n}$ and is consistent if $A$
is. Moreover, $\mathcal{E}(SAS^{-1})=S\mathcal{E}(A)$.
\end{lemma}

By Lemma \ref{lsim}, in proving some of our results, we may assume that the
reciprocal matrix has a convenient form obtained by monomial similarity.

\bigskip

For $A\in\mathcal{PC}_{n}$ and $w\in\mathbb{R}_{+}^{n}$, there is a graph
theoretic method to decide whether $w\in\mathcal{E}(A).$ Let $G(A,w)$ be the
digraph with vertex set $N$ and a directed edge $(i,j)$ if and only if
$w_{i}\geq a_{ij}w_{j}$, $1\leq i\neq j\leq n$.

\begin{theorem}
\cite{blanq2006,FJ2,FJ3}\label{blanq} Suppose that $A\in\mathcal{PC}_{n}$ and
$w\in\mathbb{R}_{+}^{n}$. The following are equivalent:

\begin{enumerate}
\item[(i)] $w\in\mathcal{E}(A)$;

\item[(ii)] $G(A,w)$ is strongly connected;

\item[(iii)] $G(A,w)$ contains a Hamiltonian cycle.
\end{enumerate}
\end{theorem}

The following observation was made in \cite{FJchap}.

\begin{remark}
\label{remspannintree}Suppose that $A\in\mathcal{PC}_{n}$. A labelled tree
(connected acyclic graph) $T$ on $n$ vertices may be viewed as a spanning tree
of the complete undirected graph with $n$ vertices. There is a unique
consistent matrix $ww^{(-T)}$ whose entries in positions given by the edges of
$T$ coincide with those of $A$ in the same positions (entries $a_{ij}$ and
$a_{ji}$ for the undirected edge $\{i,j\}$). By Theorem \ref{blanq}, vector
$w$ is efficient for $A$ because, if $\{i,j\}$ is an edge in $T$, then both
$(i,j)$ and $(j,i)$ are edges in $G(A,w),$ implying that $G(A,w)$ is strongly
connected. So, there is an element of $\mathcal{E}(A)$ associated with any
spanning tree of the complete graph. Examples are the columns of $A$
\cite{FJ1}, which correspond to the spanning trees that are labelled stars.
Another example is a vector $w$ such that $\frac{w_{\gamma_{i}}}%
{w_{\gamma_{i+1}}}=a_{\gamma_{i},\gamma_{i+1}}$ for all $i\in\{1,\ldots
,n-1\},$ in which $\gamma_{1}\cdots\gamma_{n}$ is a path with $n$ vertices in
the complete graph.
\end{remark}

\subsection{$\mathcal{E}(A)$ as a union of convex sets\label{sunion}}

For $A\in\mathcal{PC}_{n}^{0}$, the set $\mathcal{E}(A)$ is known to be
connected \cite{blanq2006,FJ4}. However, it is usually not convex
\cite{FJEJOR}. Next we describe $\mathcal{E}(A)$ as a union of a finite number
of constituent convex sets.

We identify the entries of $A=[a_{ij}]$ with the edges of the complete
directed graph on $n$ vertices, and view these interchangeably. Associated
with the cycle $\alpha=i_{1}i_{2}\cdots i_{k}i_{1}$ is the product of the
corresponding entries of $A$: $\alpha(A)=a_{i_{1}i_{2}}\cdots a_{i_{k}i_{1}}$.
If the cycle is Hamiltonian, its product in $A$ is called an \emph{H-cycle
product}. Such products are preserved by diagonal similarity on $A$. If we
reverse the direction of the edges in a cycle $\alpha$, we refer to the
alternate cycle as $\operatorname*{rev}\alpha$. Of course $\operatorname*{rev}%
\alpha(A)=\frac{1}{\alpha(A)},$ because $A$ is reciprocal. So, at most one of
$\alpha$ or $\operatorname*{rev}\alpha$ can have product $<1$.

We denote by $\Gamma(A)$ the collection of all H-cycles with product $<1$ in
$A$.\ We have $\#\Gamma(A)\leq\frac{(n-1)!}{2},$ and $\#\Gamma(A)\geq1$ since
$A$ is inconsistent.

The cycles in $\Gamma(A)$ are important in determining $\mathcal{E}(A)$ for
$A\in\mathcal{PC}_{n}^{0}$. If $\gamma\in\Gamma(A)$, with $\gamma=\gamma
_{1}\cdots\gamma_{n}\gamma_{1}$, we denote by $\mathcal{E}_{\gamma}(A)$ the
set of vectors $w=\left[
\begin{array}
[c]{ccc}%
w_{1} & \cdots & w_{n}%
\end{array}
\right]  ^{T}\in\mathbb{R}_{+}^{n}$ such that
\begin{equation}
w_{\gamma_{1}}\geq a_{\gamma_{1}\gamma_{2}}w_{\gamma_{2}}\geq a_{\gamma
_{1}\gamma_{2}}a_{\gamma_{2}\gamma_{3}}w_{\gamma_{3}}\geq\cdots\geq\left(
{\displaystyle\prod\limits_{i=1}^{n-1}}
a_{\gamma_{i}\gamma_{i+1}}\right)  w_{\gamma_{n}}\geq\gamma(A)w_{\gamma_{1}}.
\label{in}%
\end{equation}
Notice that $\mathcal{E}_{\gamma}(A)$ just depends on the entries of $A$ along
$\gamma$.

\bigskip

Recently, the following description of $\mathcal{E}(A)$ was given.

\begin{theorem}
\cite{FJ3}\label{thunion} If $A\in\mathcal{PC}_{n}^{0}$,
\[
\mathcal{E}(A)=%
{\displaystyle\bigcup\limits_{\gamma\in\Gamma_{A}}}
\mathcal{E}_{\gamma}(A).
\]

\end{theorem}

Since the solution set of the linear inequalities (\ref{in}) in $\mathbb{R}%
_{+}^{n}$ is convex, this gives $\mathcal{E}(A)$ as the union of a finite
number of convex sets. Each set has infinitely many (projectively distinct) members.

\subsection{Matrix intervals for consistent matrices built from $\mathcal{E}%
(A)$}

Let $A\in\mathcal{PC}_{n}^{0}$. Within $\gamma\in\Gamma_{A}$, for each pair
$i,j$, $i\neq j$, there is a directed path from $i$ to $j$. Denote the product
of the entries in $A$ corresponding to the edges in this path by $P_{A,\gamma
}(i,j)$. If $(i,j)$ is an edge in $\gamma$, then $P_{A,\gamma}(i,j)=a_{ij}$.
Also,%
\begin{equation}
P_{A,\gamma}(i,j)P_{A,\gamma}(j,i)=\gamma(A)<1. \label{Pijji}%
\end{equation}
We denote by $P_{A,\gamma}$ the $n$-by-$n$ matrix with $1$'s on the diagonal
and the $i,j$ entry $P_{A,\gamma}(i,j)$, $i\neq j$.

Throughout, inequalities between matrices are to be interpreted entrywise.
Recently, the following was shown.

\begin{theorem}
\cite{FJ8}\label{Thintgamma} If $A\in\mathcal{PC}_{n}^{0}$, $\gamma\in
\Gamma(A)$ and $w\in\mathcal{E}_{\gamma}(A)$, then
\begin{equation}
P_{A,\gamma}\leq W\leq P_{A,\gamma}^{(-T)}. \label{int}%
\end{equation}
Moreover, if $W=ww^{(-T)}$ satisfies (\ref{int}), then $w\in\mathcal{E}%
_{\gamma}(A)$.
\end{theorem}

Note that, because of (\ref{Pijji}), the interval in (\ref{int}) corresponding
to any $i,j$ position, $i\neq j$, is not a singleton. Moreover, each of the
inequalities in (\ref{int}) implies the other.

\bigskip

The collection of intervals (\ref{int}), with $\gamma\in\Gamma(A)$, for a
consistent matrix obtained from an efficient vector is a refinement of earlier
work \cite{FJ80} giving a single matrix interval that was the best, but not
quite as precise. This result can now be viewed as a corollary of Theorem
\ref{Thintgamma}.

Let $L_{A}$ denote the entrywise minimum (over $\gamma\in\Gamma(A)$) of
$P_{A,\gamma}$. Corresponding to (\ref{Pijji}), we also have
\[
L_{A}(i,j)L_{A}(j,i)<1.
\]

\begin{corollary}
\cite{FJ80}\label{Thintgea} If $A\in\mathcal{PC}_{n}^{0}$ and $w\in
\mathcal{E}(A)$, then%
\begin{equation}
L_{A}\leq W\leq L_{A}^{(-T)}\text{.} \label{bLA}%
\end{equation}
Moreover, $L_{A}\leq A\leq L_{A}^{(-T)}$.
\end{corollary}

However, not all consistent $W$ satisfying (\ref{bLA}) are built from
efficient vectors \cite{FJ80}.

\bigskip

Because of the next result, the matrix intervals given in Theorem
\ref{Thintgamma} and Corollary \ref{Thintgea} are the best possible.

\begin{theorem}
\cite{FJ80,FJ8}\label{rattain} If $A\in\mathcal{PC}_{n}^{0}$, $\gamma\in
\Gamma(A)$ and $1\leq i,j\leq n$, $i\neq j$, there are $w,v\in\mathcal{E}%
_{\gamma}(A)$ such that $\frac{w_{i}}{w_{j}}=P_{A,\gamma}(i,j)$ and
$\frac{v_{i}}{v_{j}}=\frac{1}{P_{A,\gamma}(j,i)}$ (the $i,j$ entry of
$P_{A,\gamma}^{(-T)}$). In particular, there are $w,v\in\mathcal{E}(A)$ such
that $\frac{w_{i}}{w_{j}}=L_{A}(i,j)$ and $\frac{v_{i}}{v_{j}}=\frac{1}%
{L_{A}(j,i)}$.
\end{theorem}

\section{Additional observations about efficient vectors\label{saddobs}}

In this section, we present several, mostly new, observations about efficient
vectors for a reciprocal matrix.

A graph characterization of when $w\in\mathbb{R}_{+}^{n}$ is efficient for a
given $A\in\mathcal{PC}_{n}$ was presented in Theorem \ref{blanq}. Next we
give an alternate (though equivalent) characterization using explicit
conditions on $A$ and $w$. If $A$ were consistent, then $w\in\mathcal{E}(A)$
if and only if $w$ is a positive multiple of any column of $A$.

\begin{proposition}
Let $A\in\mathcal{PC}_{n}^{0}$, $w\in\mathbb{R}_{+}^{n}$ and
$D=\operatorname*{diag}(w).$ Then, $w\in\mathcal{E}(A)$ if and only if
$D^{-1}AD$ contains an H-cycle with all entries $\leq1$ and at least one entry
$<1$.
\end{proposition}

\begin{proof}
By Lemma \ref{lsim}, $w\in\mathcal{E}(A)$ if and only if $\mathbf{e}%
=D^{-1}w\in\mathcal{E}(B)$, with $B=D^{-1}AD$. By Theorem \ref{thunion}, the
latter occurs if and only if there is a $\gamma\in\Gamma(B)$ such that
$\mathbf{e}\in\mathcal{E}_{\gamma}(B).$ Condition $\mathbf{e}\in
\mathcal{E}_{\gamma}(B)$ is equivalent to $1\geq b_{ij}$ for all edges $(i,j)$
in $\gamma$. Since $\gamma(B)=\gamma(A)<1$, the claim follows.

One could also argue that, by Theorem \ref{Thintgamma}, $\mathbf{e}%
\in\mathcal{E}_{\gamma}(B)$ if and only if all entries of $P_{B,\gamma}$ are
$\leq1$. All path products along $\gamma$ are $\leq1$ if and only if all
entries of $B$ along $\gamma$ are $\leq1$.
\end{proof}

\bigskip

Any vector $w\in\mathbb{R}_{+}^{n}$ is efficient for some $A\in\mathcal{PC}%
_{n}.$ For example, for a matrix $A$ whose $i$th column is $\frac{1}{w_{i}}w$
(the normalization of $w$ to have the $i$th entry equal to $1$). Since, by a
positive diagonal similarity, the $i$th column of a reciprocal matrix can be
made any arbitrary positive vector with $i$th entry $1$, it follows that, for
any $A\in\mathcal{PC}_{n}$, there is a matrix $A^{\prime}$ similar, by a
positive diagonal matrix, to $A$ such that $w\in\mathcal{E}(A^{\prime})$.

Also, any vector $w\in\mathbb{R}_{+}^{n}$ is efficient for some matrix
$A^{\prime}\in\mathcal{PC}_{n}$ obtained from $A$ by changing a certain pair
of entries in symmetric positions. This observation was made in \cite{Ro}
without reference to a known result from which it follows immediately. In
fact, by Theorem \ref{thRedeis}, $G(A,w)$ contains an H-path. Without loss of
generality, suppose that it is the path $12\cdots n.$ Let $A^{\prime}%
=[a_{ij}^{\prime}]\in\mathcal{PC}_{n}$ be such that $a_{n1}^{\prime}$
satisfies $w_{n}\geq a_{n1}^{\prime}w_{1}$ ($a_{1n}^{\prime}=\frac{1}%
{a_{n1}^{\prime}}$) and $a_{ij}^{\prime}=a_{ij}$ otherwise. Then the digraph
$G(A^{\prime},w)$ contains the H-cycle $12\cdots n1$. By Theorem \ref{blanq},
$w$ is efficient for $A^{\prime}$. However, if $a_{n1}$ is very large compared
with $\frac{w_{n}}{w_{1}},$ the required change in $A$ is large, unless there
are other H-paths in $G(A,w)$ that lead to lesser changes.

\begin{example}
Let
\[
A=\left[
\begin{array}
[c]{ccc}%
1 & 4 & 4\\
\frac{1}{4} & 1 & 4\\
\frac{1}{4} & \frac{1}{4} & 1
\end{array}
\right]  .
\]
Then $G(A,\mathbf{e})$ has an edge $i\rightarrow j$ if and only if $i>j.$
Thus, $321$ is the only H-path in $G(A,\mathbf{e})$. This implies that any
reciprocal matrix differing from $A$ in just one pair of entries and for which
$\mathbf{e}$ is efficient should be obtained by modifying the entries in
positions $1,3$ and $3,1.$ More precisely, the pair of entries $(4,\frac{1}%
{4})$ corresponding to these positions should be replaced by a pair
$(x,\frac{1}{x}),$ with $x\leq1$, so that $1\rightarrow3$ is an edge in
$G(A^{\prime},\mathbf{e})$, in which $A^{\prime}$ is the modified matrix.
\end{example}

\bigskip

Theorem \ref{Thintgamma} and Corollary \ref{Thintgea} shed light on the
question about the possible orders among the efficient vectors for
$A\in\mathcal{PC}_{n}^{0}$. Given $A\in\mathcal{PC}_{n}^{0}$ and $\gamma
\in\Gamma(A)$, necessary and sufficient conditions for all vectors in the set
$\mathcal{E}_{\gamma}(A)$ to have the same order (in which case $A$ is said to
be uniformly ordered for $\gamma$) are given in \cite{FJ8,FJRR}$.$ Necessary,
and also sufficient, conditions for the existence of a given order among the
efficient vectors for $A$ are discussed in \cite{FJRR}. An interesting
observation that generalizes Theorem 13 in \cite{FJRR} is the following.

\begin{theorem}
\label{Torder}Let $A\in\mathcal{PC}_{n}^{0}$ and $\gamma\in\Gamma(A).$ All
strict orders occur in $\mathcal{E}_{\gamma}(A)$ if and only if all entries in
$A$ along $\gamma$ are $<1$.
\end{theorem}

\begin{proof}
Suppose that all entries in $A$ along $\gamma$ are $<1$. Then, all
off-diagonal entries in $P_{A,\gamma}$ are $<1$ implying that $P_{A,\gamma
}<J<P_{A,\gamma}^{(-T)}$ (except for the diagonal entries). Then,
$P_{A,\gamma}\leq ww^{(-T)}\leq P_{A,\gamma}^{(-T)}$ for any sufficiently
small perturbation $w$ of $\mathbf{e}$ (which may have any arbitrary strict
order). By Theorem \ref{Thintgamma}, $w\in\mathcal{E}_{\gamma}(A)$ (this
conclusion may not hold if $\mathcal{E}(A)$ is considered instead of
$\mathcal{E}_{\gamma}(A)$). So any strict order occurs in $\mathcal{E}%
_{\gamma}(A)$ and, thus, in $\mathcal{E}(A).$

Suppose that any strict order occurs in $\mathcal{E}_{\gamma}(A)$. Then, for
each $i,j\in N$, $i\neq j,$ there are $w,v\in\mathcal{E}_{\gamma}(A)$ such
that $\frac{w_{i}}{w_{j}}>1$ and $\frac{v_{i}}{v_{j}}<1.$ By Theorem
\ref{Thintgamma}, this implies that $P_{A,\gamma}(i,j)<1$ and $P_{A,\gamma
}^{(-T)}(i,j)>1$. Thus, $P_{A,\gamma}$ has all off-diagonal entries $<1$
implying that all entries in $A$ along $\gamma$ are $<1$.
\end{proof}

\bigskip

Since all entries in $A$ along an H-cycle $\gamma$ being $<1$ implies that
$\gamma\in\Gamma(A),$ we have the following consequence of Theorem
\ref{Torder}.

\begin{corollary}
Let $A\in\mathcal{PC}_{n}^{0}$. If there is an H-cycle along which all entries
in $A$ are $<1$, then all strict orders occur in $\mathcal{E}(A).$
\end{corollary}

\bigskip

An $n$-by-$n$ matrix $A$ is said to be a \emph{partial reciprocal matrix}
(PRM) if some of its entries are known and others are unspecified and, in
addition, if the entry $a_{ij}$ is specified then so is $a_{ji}$ and
$a_{ij}=\frac{1}{a_{ji}}.$ A PRM admits a reciprocal completion. Associate
with a PRM $A$ the undirected graph with vertex set $N$ and an edge $\{i,j\}$
if and only if $a_{ij}$ (and thus $a_{ji}$) is specified. We call it the
\emph{specified-entry graph} of $A$ and denote it by $S_{A}$.

\begin{theorem}
\label{Tinc}Let $A$ be an $n$-by-$n$ PRM. There is a $w\in\mathbb{R}_{+}^{n}$
efficient for all completions $B$ of $A$ if and only if $S_{A}$ is connected.
\end{theorem}

\begin{proof}
($\Leftarrow$) Suppose that $S_{A}$ is connected. Let $T$ be a spanning tree
of $S_{A}$. Let $w$ be such that the consistent matrix $ww^{(-T)}$ coincides
with $A$ in the entries given by $T$. Then $G(B,w)$ is strongly connected for
any completion $B$ of $A$.

($\Rightarrow$) Let $B$ be a completion of $A$ and $w\in\mathbb{R}_{+}^{n}$ be
efficient for $B$. Suppose that $S_{A}$ is not connected. We will see that
there is another completion $C$ of $A$ such that $w$ is not efficient for $C$.
Let $\{S_{1},S_{2}\}$ be a partition of $N$ such that there is no edge in
$S_{A}$ between a vertex in $S_{1}$ and a vertex in $S_{2}.$ For all $i\in
S_{1}$ and $j\in S_{2}$, let $c_{ij}>0$ be such that $\frac{w_{i}}{w_{j}%
}>c_{ij}$ and let $c_{ji}=\frac{1}{c_{ij}}.$ Let $c_{ij}=a_{ij}$ if the $i,j$
entry of $A$ is specified. Choose the remaining reciprocal pairs of
unspecified entries arbitrarily to get a reciprocal completion $C$ of $A.$
Then $S_{2}$ is a sink set of $G(C,w)$ and thus $w$ is not efficient for $C.$
\end{proof}

\section{Efficient vectors built from efficient vectors for principal
submatrices\label{Seffprincsubm}}

In \cite{FJ1} it was shown that, for $S\subset N,$ any efficient vector for
$A[S]$ can be extended to a vector in $\mathcal{E}(A)$, that is, for any
$u\in\mathcal{E}(A[S]),$ there is a $w\in\mathcal{E}(A)$ such that $w[S]=u.$
The possible extensions were characterized. An example is the following. Any
weighted geometric mean of the columns of $A[S]$ is efficient for $A[S]$
\cite{FJ2}. Also, the weighted geometric mean of the columns of $A$ with the
same weights for the columns indexed by $S\ $(and $0$ weights for the
remaining columns) is efficient for $A$. The latter vector is an extension of
the former.

\bigskip Generalizing the result of \cite{FJ1}, we next show that, for any
partition $\{S_{1},\ldots,S_{p}\}$ of $N$ and any vectors $u^{[i]}%
\in\mathcal{E}(A[S_{i}])$, $i=1,\ldots p,$ there is a $w\in\mathcal{E}(A)$
such that $w[S_{i}]$ is proportional to $u^{[i]}$, that is, they determine the
same consistent matrix. We describe, inductively, how these vectors are constructed.

\begin{theorem}
\label{Teffprinc}Let $A\in\mathcal{PC}_{n}$. Let $S_{1}=\{1,\ldots,k\}$ and
$S_{2}=\{k+1,\ldots,n\}$, so that $\{S_{1},S_{2}\}$ is a partition of $N.$
Let
\[
u^{[1]}=\left[
\begin{array}
[c]{ccc}%
u_{1} & \cdots & u_{k}%
\end{array}
\right]  ^{T}\in\mathcal{E}(A[S_{1}])\qquad\text{and}\qquad u^{[2]}=\left[
\begin{array}
[c]{ccc}%
u_{k+1} & \cdots & u_{n}%
\end{array}
\right]  ^{T}\in\mathcal{E}(A[S_{2}]).
\]
Let
\[
F=\left\{  a_{ij}\frac{u_{j}}{u_{i}}:i\in S_{1}\text{ and }j\in S_{2}\right\}
,
\]
$c_{\max}=\max F$ and $c_{\min}=\min F.$ Then, the following are equivalent:

\begin{enumerate}
\item[(i)] $w\in\mathcal{E}(A)$, $W[S_{1}]=u^{[1]}\left(  u^{[1]}\right)
^{(-T)}$, and $W[S_{2}]=u^{[2]}\left(  u^{[2]}\right)  ^{(-T)}$;

\item[(ii)] up to a positive factor,
\begin{equation}
w=\left[
\begin{array}
[c]{c}%
cu^{[1]}\\
u^{[2]}%
\end{array}
\right]  , \label{wef}%
\end{equation}
for some $c\in\lbrack c_{\min},c_{\max}]$.
\end{enumerate}
\end{theorem}

\begin{proof}
((i)$\Rightarrow$(ii)) The last two conditions in (i) imply that $w[S_{1}]$ is
a multiple of $u^{[1]}$ and $w[S_{2}]$ is a multiple of $u^{[2]}$. Then, to
describe $w$ up to a positive factor, we may assume that $w[S_{2}]=u^{[2]}.$
Let $c>0$ be such that $w[S_{1}]=cu^{[1]}$. Since, by Theorem \ref{blanq},
$G(A,w)$ is strongly connected, there are $i_{0},i_{0}^{\prime}\in S_{1}$ and
$j_{0},j_{0}^{\prime}\in S_{2}$ such that $(i_{0},j_{0})$ and $(j_{0}^{\prime
},i_{0}^{\prime})$ are edges in $G(A,w)$. Then,%
\[
\frac{cu_{i_{0}}}{u_{j_{0}}}=\frac{w_{i_{0}}}{w_{j_{0}}}\geq a_{i_{0}j_{0}}%
\]
and%
\[
\frac{cu_{i_{0}^{\prime}}}{u_{j_{0}^{\prime}}}=\frac{w_{i_{0}^{\prime}}%
}{w_{j_{0}^{\prime}}}\leq a_{i_{0}^{\prime}j_{0}^{\prime}}.
\]
This implies that
\[
a_{i_{0}j_{0}}\frac{u_{j_{0}}}{u_{i_{0}}}\leq c\leq a_{i_{0}^{\prime}%
j_{0}^{\prime}}\frac{u_{j_{0}^{\prime}}}{u_{i_{0}^{\prime}}}.
\]
Since $c_{\min}\leq a_{i_{0}j_{0}}\frac{u_{j_{0}}}{u_{i_{0}}}$ and $c_{\max
}\geq a_{i_{0}^{\prime}j_{0}^{\prime}}\frac{u_{j_{0}^{\prime}}}{u_{i_{0}%
^{\prime}}}$, the implication follows.

((ii)$\Rightarrow$(i)) Let $w$ be as in (\ref{wef}) for some $c\in\lbrack
c_{\min},c_{\max}]$. Clearly, the last two conditions in (i) hold. Let
$i_{0},i_{0}^{\prime}\in S_{1}$ and $j_{0},j_{0}^{\prime}\in S_{2}$ be such
that $a_{i_{0}j_{0}}\frac{u_{j_{0}}}{u_{i_{0}}}\leq c$ and $a_{i_{0}^{\prime
}j_{0}^{\prime}}\frac{u_{j_{0}^{\prime}}}{u_{i_{0}^{\prime}}}\geq c.$ Then
\[
\frac{w_{i_{0}}}{w_{j_{0}}}=\frac{cu_{i_{0}}}{u_{j_{0}}}\geq a_{i_{0}j_{0}%
}\frac{u_{j_{0}}}{u_{i_{0}}}\frac{u_{i_{0}}}{u_{j_{0}}}=a_{i_{0}j_{0}}%
\]
and%
\[
\frac{w_{i_{0}^{\prime}}}{w_{j_{0}^{\prime}}}=\frac{cu_{i_{0}^{\prime}}%
}{u_{j_{0}^{\prime}}}\leq a_{i_{0}^{\prime}j_{0}^{\prime}}\frac{u_{j_{0}%
^{\prime}}}{u_{i_{0}^{\prime}}}\frac{u_{i_{0}^{\prime}}}{u_{j_{0}^{\prime}}%
}=a_{i_{0}^{\prime}j_{0}^{\prime}}.
\]
Thus, $(i_{0},j_{0})$ and $(j_{0}^{\prime},i_{0}^{\prime})$ are edges in
$G(A,w)$. On the other hand, $G(A[S_{1}],cu^{[1]})$ and $G(A[S_{2}],u^{[2]})$
are induced subgraphs of $G(A,w)$ and are strongly connected by Theorem
\ref{blanq}. Note that $G(A[S_{1}],u^{[1]})$ and $G(A[S_{1}],cu^{[1]})$
coincide. Thus, $G(A,w)$ is strongly connected, implying, by Theorem
\ref{blanq}, that $w\in\mathcal{E}(A).$
\end{proof}

\bigskip

The following extension to a more refined partition of $N$ follows from
Theorem \ref{Teffprinc} using induction on the number of sets in the partition
(and, if necessary, considering a permutation similarity of the original
reciprocal matrix).

\begin{theorem}
\label{thprincipalblocks}Let $A\in\mathcal{PC}_{n}$. Let $\{S_{1},\ldots
,S_{p}\}$ be a partition of $N.$ Let $w^{[i]}\in\mathcal{E}(A[S_{i}])$, for
$i=1,\ldots,p.$ Then, there is a $w\in\mathcal{E}(A)$ such that $W[S_{i}%
]=w^{[i]}(w^{[i]})^{(-T)}$, that is, $w[S_{i}]$ is a positive multiple of
$w^{[i]}$, for $i=1,\ldots,p$.
\end{theorem}

The particular form of the vectors $w$ in Theorem \ref{thprincipalblocks}
follows from repeated use of the form in Theorem \ref{Teffprinc}.

\bigskip

When the comparisons in $A$ within disjoint sets of alternatives are
consistent, one may want to choose an efficient vector that preserves them.
Theorem \ref{thprincipalblocks} guarantees that this is possible.

\bigskip We note that Theorem \ref{thprincipalblocks} can also be obtained
using the definition of efficiency directly. If we start with a vector $w$
whose subvectors are efficient for the corresponding principal submatrices,
then any efficient vector dominating $w$ (whose existence is claimed, for
example, in \cite{BoFU}) should have the corresponding subvectors dominating
those in $w$ and, thus, they should be multiples. However, our proof of
Theorems \ref{Teffprinc} and \ref{thprincipalblocks} gives an explicit method
to obtain such a $w.$

\section{Efficient vectors preserving orders among sets of alternatives
\label{Spreservord}}

We may be interested in efficient vectors $w$ for $A\in\mathcal{PC}_{n}^{0}$
that satisfy prescribed order relations among certain alternatives. More
precisely, given a set $E\subset N\times N$ of pairs with distinct entries, we
may require that, for any $(i,j)\in E$, $w_{i}\geq w_{j}.$ Of course "order
consistency" should occur. However, this is not sufficient as, in general, not
any order can occur among the efficient vectors (see Theorem \ref{Torder}). In
view of Theorem \ref{thunion}, the set of efficient vectors satisfying the
imposed inequalities is a (possibly empty) union of convex sets. The order
restrictions may be internal (depending on $A$) or externally imposed. Here we
focus on the former case.

\bigskip

The particular case of an efficient vector preserving the orders induced by
$A$ along a spanning tree is especially simple. In fact, in view of Remark
\ref{remspannintree}, more holds. There is a positive vector $w$ such that $A$
and $W$ coincide along the tree. We consider here an extension of this order
statement to sets of alternatives. We first give a more general result.

\bigskip

We denote by $T_{p}$ a spanning tree of the complete labelled graph with $p$
vertices. For convenience, we allow $p$ to be $1$, in which case $T_{p}$ has
just one vertex.

\begin{theorem}
\label{thstep}Let $\{S_{1},\ldots,S_{p}\}$ be a partition of $N.$ Let
$A\in\mathcal{PC}_{n}$ and let $w^{[i]}\in\mathcal{E}(A[S_{i}])$,
$i=1,\ldots,p.$ Then, there is a $w\in\mathcal{E}(A)$ such that $W[S_{i}%
]=w^{[i]}(w^{[i]})^{(-T)}$, $i=1,\ldots,p$. Such a $w$ may be chosen so that,
for any edge $\{i,j\}$ in a particular $T_{p}$, any one of the following
conditions holds:

\begin{enumerate}
\item[(i)] $\frac{w_{q_{i}}}{w_{q_{j}}}=a_{q_{i}q_{j}}$, for a given $q_{i}\in
S_{i}$ and a given $q_{j}\in S_{j}$;

\item[(ii)] $W[S_{i},S_{j}]\geq A[S_{i},S_{j}]$;

\item[(iii)] $W[S_{i},S_{j}]\leq A[S_{i},S_{j}]$.
\end{enumerate}
\end{theorem}

\begin{proof}
The proof is by induction on $p.$ When $p=1$, $w^{[1]}\in\mathcal{E}(A)$, so
that we may take $w$ to be any positive multiple of this vector (which will
lie in $\mathcal{E}(A)$, as required). Suppose that $p>1.$ Relabelling the
sets of the partition, we may assume, without loss of generality, that $p$ is
a pendent vertex and $\{p-1,p\}$ is an edge in $T_{p}$. Let $S^{0}=S_{1}%
\cup\cdots\cup S_{p-1}$ and $T_{p-1}$ be the tree obtained from $T_{p}$ by
deleting vertex $p$. By a possible permutation similarity on $A$, we may
assume in addition that $S^{0}=\{1,\ldots,k\}$ and $S_{p}=\{k+1,\ldots,n\},$
with $1\leq k<n$.

By the induction hypothesis, there is $u=\left[
\begin{array}
[c]{ccc}%
u_{1} & \cdots & u_{k}%
\end{array}
\right]  ^{T}\in\mathcal{E}(A[S^{0}])$ such that $\left(  uu^{(-T)}\right)
[S_{i}]=w^{[i]}(w^{[i]})^{(-T)},$ $i=1,\ldots,p-1$,$\ $ and, in addition, for
any edge $\{i,j\}$ in $T_{p-1}$, any one of conditions (i), (ii), or (iii) may occur.

Let $w^{[p]}=\left[
\begin{array}
[c]{ccc}%
u_{k+1} & \cdots & u_{n}%
\end{array}
\right]  ^{T}$. For $q_{p-1}\in S_{p-1}$ and $q_{p}\in S_{p},$ let
\[
c=a_{q_{p-1}q_{p}}\frac{u_{q_{p}}}{u_{q_{p-1}}}\text{.}%
\]
Since
\[
\min\left\{  a_{ij}\frac{u_{j}}{u_{i}}:i\in S^{0}\text{ and }j\in
S_{p}\right\}  \leq c\leq\max\left\{  a_{ij}\frac{u_{j}}{u_{i}}:i\in
S^{0}\text{ and }j\in S_{p}\right\}  ,
\]
by Theorem \ref{Teffprinc},
\[
w=\left[
\begin{array}
[c]{c}%
cu\\
w^{[p]}%
\end{array}
\right]  \in\mathcal{E}(A).
\]
Note that $(cu)(cu)^{(-T)}=uu^{(-T)}.$

We have%
\[
\frac{w_{q_{p-1}}}{w_{q_{p}}}=\frac{cu_{q_{p-1}}}{u_{q_{p}}}=a_{q_{p-1}q_{p}%
}.
\]
If $q_{p-1}\in S_{p-1}$ and $q_{p}\in S_{p}$ are such that%
\[
c=a_{q_{p-1}q_{p}}\frac{u_{q_{p}}}{u_{q_{p-1}}}=\max\left\{  a_{ij}\frac
{u_{j}}{u_{i}}:i\in S_{p-1}\text{ and }j\in S_{p}\right\}  ,
\]
then, for $i\in S_{p-1}$ and $j\in S_{p}$,
\[
\frac{w_{i}}{w_{j}}=\frac{cu_{i}}{u_{j}}\geq a_{ij}\frac{u_{j}}{u_{i}}%
\frac{u_{i}}{u_{j}}=a_{ij},
\]
implying $W[S_{i},S_{j}]\geq A[S_{i},S_{j}]$.

Similarly, if $q_{p-1}\in S_{p-1}$ and $q_{p}\in S_{p}$ are such that%
\[
c=a_{q_{p-1}q_{p}}\frac{u_{q_{p}}}{u_{q_{p-1}}}=\min\left\{  a_{ij}\frac
{u_{j}}{u_{i}}:i\in S_{p-1}\text{ and }j\in S_{p}\right\}  ,
\]
then, for $i\in S_{p-1}$ and $j\in S_{p},$
\[
\frac{w_{i}}{w_{j}}=\frac{cu_{i}}{u_{j}}\leq a_{ij}\frac{u_{j}}{u_{i}}%
\frac{u_{i}}{u_{j}}=a_{ij},
\]
implying $W[S_{i},S_{j}]\leq A[S_{i},S_{j}]$.
\end{proof}

\bigskip

For the spanning tree $T_{p}$ on $p$ vertices, we assign a weight $t_{ij}%
\in\{-1,1\}$ to each edge $\{i,j\}$. We say that $A\in\mathcal{PC}_{n}$
is$\ \{S_{1},\ldots,S_{p}\}$\emph{-order unambiguous for the weighted tree
}$T_{p}$ ($\{S_{1},\ldots,S_{p}\}$-OU for $T_{p}$) if, for any edge $\{i,j\}$
in $T_{p}$ in which $i<j$, we have $A[S_{i},S_{j}]\geq J$ if $t_{ij}=1$ and
$A[S_{i},S_{j}]\leq J$ if $t_{ij}=-1$. The former (latter) inequality means
that any alternative in $S_{i}$ is weakly ranked above (below) any alternative
in $S_{j}$, by the direct comparisons in $A$. When this terminology is used,
it is implicit that $\{S_{1},\ldots,S_{p}\}$ is a partition of $N$ and $T_{p}$
is a weighted tree on $p$ vertices.

\bigskip

The next corollary is an immediate consequence of Theorem \ref{thstep}.

\begin{corollary}
\label{corder}Let $\{S_{1},\ldots,S_{p}\}$ be a partition of $N.$ Let
$A\in\mathcal{PC}_{n}$ be $\{S_{1},\ldots,S_{p}\}$-OU for $T_{p}.$ Let
$w^{[i]}\in\mathcal{E}(A[S_{i}])$, $i=1,\ldots,p.$ Then, there is a
$w\in\mathcal{E}(A)$ such that $W[S_{i}]=w^{[i]}(w^{[i]})^{(-T)}$,
$i=1,\ldots,p$, and, in addition, $W$ is $\{S_{1},\ldots,S_{p}\}$-OU for
$T_{p}$.
\end{corollary}

\bigskip

We illustrate the construction of the vector $w$ in Corollary \ref{corder} in
the case $T_{p}$ is a path.

\begin{example}
Let%
\[
A=\left[
\begin{tabular}
[c]{cc|ccc|cc}%
$1$ & $2$ & $2$ & $3$ & $\frac{3}{2}$ & $2$ & $\frac{1}{2}$\\
$\frac{1}{2}$ & $1$ & $\frac{4}{3}$ & $2$ & $\frac{3}{2}$ & $\frac{1}{2}$ &
$2$\\\hline
$\frac{1}{2}$ & $\frac{3}{4}$ & $1$ & $2$ & $\frac{1}{3}$ & $\frac{1}{2}$ &
$\frac{1}{3}$\\
$\frac{1}{3}$ & $\frac{1}{2}$ & $\frac{1}{2}$ & $1$ & $3$ & $\frac{1}{2}$ &
$\frac{1}{2}$\\
$\frac{2}{3}$ & $\frac{2}{3}$ & $3$ & $\frac{1}{3}$ & $1$ & $\frac{1}{3}$ &
$\frac{1}{2}$\\\hline
$\frac{1}{2}$ & $2$ & $2$ & $2$ & $3$ & $1$ & $\frac{1}{4}$\\
$2$ & $\frac{1}{2}$ & $3$ & $2$ & $2$ & $4$ & $1$%
\end{tabular}
\right]  .
\]
Let $S_{1}=\{1,2\},$ $S_{2}=\{3,4,5\}$ and $S_{3}=\{6,7\}.$ Matrix $A$ is
$\{S_{1},S_{2},S_{3}\}$-OU for the path $123$ with weights $t_{12}=1$ and
$t_{23}=-1.$ Let $w^{[i]}\ $ be an efficient vector for $A[S_{i}],$ $i=1,2,3.$
We take $w^{[i]}$ the last column of $A[S_{i}]$. Write $w^{[1]}=\left[
\begin{array}
[c]{cc}%
w_{1} & w_{2}%
\end{array}
\right]  ^{T},$ $w^{[2]}=\left[
\begin{array}
[c]{ccc}%
w_{3} & w_{4} & w_{5}%
\end{array}
\right]  ^{T}$ and $w^{[3]}=\left[
\begin{array}
[c]{cc}%
w_{6} & w_{7}%
\end{array}
\right]  ^{T}.$ Let
\begin{align*}
c_{1}  &  =\max\left\{  a_{ij}\frac{w_{j}}{w_{i}}:i\in S_{1}\text{ and }j\in
S_{2}\right\} \\
&  =\max\left\{  a_{13}\frac{w_{3}}{w_{1}},a_{14}\frac{w_{4}}{w_{1}}%
,a_{15}\frac{w_{5}}{w_{1}},a_{23}\frac{w_{3}}{w_{2}},a_{24}\frac{w_{4}}{w_{2}%
},a_{25}\frac{w_{5}}{w_{2}}\right\} \\
&  =\max\left\{  2\frac{\frac{1}{3}}{2},3\frac{3}{2},\frac{3}{2}\frac{1}%
{2},\frac{4}{3}\frac{\frac{1}{3}}{1},2\frac{3}{1},\frac{3}{2}\frac{1}%
{1}\right\}  =6.
\end{align*}
Let
\[
u=\left[
\begin{array}
[c]{c}%
c_{1}w^{[1]}\\
w^{[2]}%
\end{array}
\right]  =\left[
\begin{array}
[c]{ccccc}%
12 & 6 & \frac{1}{3} & 3 & 1
\end{array}
\right]  ^{T}\in\mathcal{E}(A[S_{1}\cup S_{2}]).
\]
Let
\begin{align*}
c_{2}  &  =\min\left\{  a_{ij}\frac{w_{j}}{w_{i}}:i\in S_{1}\cup S_{2}\text{
and }j\in S_{3}\right\} \\
&  =\min\left\{  a_{16}\frac{w_{6}}{w_{1}},a_{17}\frac{w_{7}}{w_{1}}%
,a_{26}\frac{w_{6}}{w_{2}},a_{27}\frac{w_{7}}{w_{2}},a_{36}\frac{w_{6}}{w_{3}%
},a_{37}\frac{w_{7}}{w_{3}},a_{46}\frac{w_{6}}{w_{4}},a_{47}\frac{w_{7}}%
{w_{4}},a_{56}\frac{w_{6}}{w_{5}},a_{57}\frac{w_{7}}{w_{5}}\right\} \\
&  =\min\left\{  2\frac{\frac{1}{4}}{12},\frac{1}{2}\frac{1}{12},\frac{1}%
{2}\frac{\frac{1}{4}}{6},2\frac{1}{6},\frac{1}{2}\frac{\frac{1}{4}}{\frac
{1}{3}},\frac{1}{3}\frac{1}{\frac{1}{3}},\frac{1}{2}\frac{\frac{1}{4}}%
{3},\frac{1}{2}\frac{1}{3},\frac{1}{3}\frac{\frac{1}{4}}{1},\frac{1}{2}%
\frac{1}{1}\right\}  =\frac{1}{48}.
\end{align*}
Then,
\[
w=\left[
\begin{array}
[c]{c}%
c_{2}u\\
w^{[3]}%
\end{array}
\right]  =\left[
\begin{array}
[c]{ccccccc}%
\frac{1}{4} & \frac{1}{8} & \frac{1}{144} & \frac{1}{16} & \frac{1}{48} &
\frac{1}{4} & 1
\end{array}
\right]  ^{T}\in\mathcal{E}(A).
\]
We have%
\[
W=\left[
\begin{tabular}
[c]{cc|ccc|cc}%
$1$ & $2$ & $36$ & $4$ & $12$ & $1$ & $\frac{1}{4}$\\
$\frac{1}{2}$ & $1$ & $18$ & $2$ & $6$ & $\frac{1}{2}$ & $\frac{1}{8}$\\\hline
$\frac{1}{36}$ & $\frac{1}{18}$ & $1$ & $\frac{1}{9}$ & $\frac{1}{3}$ &
$\frac{1}{36}$ & $\frac{1}{144}$\\
$\frac{1}{4}$ & $\frac{1}{2}$ & $9$ & $1$ & $3$ & $\frac{1}{4}$ & $\frac
{1}{16}$\\
$\frac{1}{12}$ & $\frac{1}{6}$ & $3$ & $\frac{1}{3}$ & $1$ & $\frac{1}{12}$ &
$\frac{1}{48}$\\\hline
$1$ & $2$ & $36$ & $4$ & $12$ & $1$ & $\frac{1}{4}$\\
$4$ & $8$ & $144$ & $16$ & $48$ & $4$ & $1$%
\end{tabular}
\ \right]  ,
\]
which is $\{S_{1},S_{2},S_{3}\}$-OU for the path $123$ with weights $t_{12}=1$
and $t_{23}=-1.$ Observe that (projectively) different choices of the
efficient vectors for the principal submatrices give different matrices $W$.
\end{example}

\bigskip

The vector $w$ in Corollary \ref{corder} verifies the following claim.

\begin{corollary}
\label{cRR}Let $\{S_{1},\ldots,S_{p}\}$ be a partition of $N.$ Let
$A\in\mathcal{PC}_{n}$ be $\{S_{1},\ldots,S_{p}\}$-OU for $T_{p}$. Then, there
is a $w\in\mathcal{E}(A)$ satisfying $w[S_{i}]\in\mathcal{E}(A[S_{i}]),$
$i=1,\ldots,p,$ and such that, for all edge $\{i,j\}$ in $T_{p}$ and all
$q_{i}\in S_{i}$ and $q_{j}\in S_{j}$, there is no rank reversal in $A$ at
position $q_{i},q_{j}$.
\end{corollary}

For any weights in $T_{p}$, there is a consistent matrix $W=ww^{(-T)}$ that is
$\{S_{1},\ldots,S_{p}\}$-OU for $T_{p}$. However, in general, the order of
such a $w$ is not determined by the weights. For example, if $W$ is
$\{S_{1},S_{2},S_{3},S_{4}\}$-OU for the path $1234$ with corresponding
weights $\ 1,-1,-1$, that is, $S_{1}\succsim S_{2}$, $S_{3}\succsim S_{2},$
$S_{4}\succsim S_{3}$, we may have $S_{1}\succsim S_{4}\succsim S_{3}\succsim
S_{2},$ $S_{4}\succsim S_{3}\succsim S_{1}\succsim S_{2}$ or $S_{4}\succsim
S_{1}\succsim S_{3}\succsim S_{2}.\ $Here, $S_{i}\succsim S_{j}$ means that
$W[S_{i},S_{j}]\geq J$. However, if the tree is a path we have the following.

\begin{corollary}
\label{cortrans}Let $\{S_{1},\ldots,S_{p}\}$ be a partition of $N$ and let
$T_{p}$ be the weighted path $12\cdots p$ with all weights $1.$ Let
$A\in\mathcal{PC}_{n}$ be $\{S_{1},\ldots,S_{p}\}$-OU for $T_{p}$. Let
$w^{[i]}\in\mathcal{E}(A[S_{i}])$, $i=1,\ldots,p.$ Then, there is a
$w\in\mathcal{E}(A)$ such that $W[S_{i}]=w^{[i]}(w^{[i]})^{(-T)},$
$i=1,\ldots,p$, and $W[S_{j},S_{j+1}]\geq J$ for all $j\in\{1,\ldots,p-1\}$.
In this case, also $W[S_{j},S_{j+1}\cup\cdots\cup S_{p}]\geq J$ for all
$j\in\{1,\ldots,p-1\}.$
\end{corollary}

\begin{proof}
The first claim is an immediate consequence of Corollary \ref{corder}. The
second claim follows since $W$ is consistent.
\end{proof}

\bigskip

Theorem \ref{Teffprinc} and Corollary \ref{corder} can be combined to study a
more general pattern of preserved orders. Let $A\in\mathcal{PC}_{n}$ and let
$\{S_{1},\ldots,S_{p}\}$ be a partition of $N$. For $1\leq r<p$, define
\[
S^{(1)}=S_{1}\cup\cdots\cup S_{r}\quad\text{and}\quad S^{(2)}=S_{r+1}%
\cup\cdots\cup S_{p}.
\]
Suppose that $A[S^{(1)}]$ is $\{S_{1},\ldots,S_{r}\}$-OU for $T_{r}$, and that
$A[S^{(2)}]$ is $\{S_{r+1},\ldots,S_{p}\}$-OU for $T_{p-r}$, while no
particular structure is assumed for the block $A[S^{(1)},S^{(2)}]$. Then, by
Corollary \ref{corder}, there exists a vector $w\in\mathbb{R}_{+}^{n}$ such
that $w[S_{i}]\in\mathcal{E}(A[S_{i}])$ for $i=1,\ldots,p$, $w[S^{(1)}]$ is
$\{S_{1},\ldots,S_{r}\}$-OU for $T_{r}$, and $w[S^{(2)}]$ is $\{S_{r+1}%
,\ldots,S_{p}\}$-OU for $T_{p-r}$. By Theorem \ref{Teffprinc}, $w$ can be
chosen in $\mathcal{E}(A)$. If $A[S^{(1)},S^{(2)}]\geq J$ or $A[S^{(1)}%
,S^{(2)}]\leq J$, applying Corollary \ref{corder} with $w^{[1]}=w[S^{(1)}]$
and $w^{[2]}=w[S^{(2)}]$ gives a $w\in\mathcal{E}(A)$ with the additional
property that $W[S^{(1)},S^{(2)}]$ is $\geq J$ or $\leq J$, as in $A$. 

\bigskip

A result corresponding to Corollary \ref{corder} may hold for some partial
reciprocal matrices. Let $\{S_{1},\ldots,S_{p}\}$ be a partition of $N$ and
$T_{p}$ be a spanning tree of the complete graph with $p$ vertices. It may
happen that all pairwise comparisons between alternatives in $S_{i}$ and
alternatives in $S_{j}$ are known for every edge $\{i,j\}$ in $T_{p}$, whereas
some other comparisons remain unspecified. This leads to a PRM matrix $A$
whose blocks $A[S_{i},S_{j}]$ are fully specified. If, for each edge
$\{i,j\}$, either $A[S_{i},S_{j}]\geq J$ or $A[S_{i},S_{j}]\leq J$, then
Theorem \ref{thstep} guarantees that, for every reciprocal completion of $A$,
there exists an efficient vector $w$ such that $W$ is $\{S_{1},\ldots,S_{p}%
\}$-OU for $T_{p}$ with respect to the weights induced by $A$. In other words,
$w$ preserves the preference ordering between any alternative in $S_{i}$ and
any alternative in $S_{j}$ implied by the direct comparisons. In general, the
efficient vector $w$ depends on the particular completion. However, the
following holds.

\begin{corollary}
Let $\{S_{1},\ldots,S_{p}\}$ be a partition of $N$ and $T_{p}$ be a given
weighted spanning tree (as above). Let $A$ be an $n$-by-$n$ partial reciprocal
matrix. Suppose that, for any edge $\{i,j\}$ in $T_{p}$, $A[S_{i},S_{j}]$ is
specified and either $A[S_{i},S_{j}]\geq J$ or $A[S_{i},S_{j}]\leq J$,
according to the corresponding weight. If the specified-entry graph associated
with each principal partial submatrix $A[S_{i}]$ is connected, then there
exists a vector $w$ efficient for all reciprocal completions of $A$ and such
that $W$ is $\{S_{1},\ldots,S_{p}\}$-OU for $T_{p}.$
\end{corollary}

\begin{proof}
If the specified-entry graph associated with each principal submatrix
$A[S_{i}]$ is connected, by Theorem \ref{Tinc}, for each $i=1,\ldots,p$, there
is a vector $w^{[i]}$ efficient for all reciprocal completions of $A[S_{i}]$.
Then, the construction in the proof of Theorem \ref{thstep} of a vector $w$
efficient for a completion of $A$ and such that $W$ is $\{S_{1},\ldots
,S_{p}\}$-OU for $T_{p}$ just depends on the vectors $w^{[i]}$ and the
specified blocks $A[S_{i},S_{j}]$ corresponding to the edges $\{i,j\}$ of
$T_{P}.$
\end{proof}

\section{Block order transitive reciprocal matrices\label{SOT}}

A matrix $A\in\mathcal{PC}_{n}$ is said to be \emph{order transitive} if
$a_{ij}\geq1$ and $a_{jk}\geq1$ imply $a_{ik}\geq1$, for all $1\leq i,j,k\leq
n$. Order transitivity is invariant under permutation similarity.

Taking into account Theorem \ref{Prop1}, it follows that, if $A$ is order
transitive, then it is permutationally similar to a matrix with all entries
above the diagonal $\geq1$. The converse is also true if all off-diagonal
entries of $A$ are distinct from $1$ \cite{FJRR}. In fact, suppose, without
loss of generality, that $A$ has all entries above the diagonal $>1$. Then,
$a_{ij}>1$ and $a_{jk}>1$ imply $j>i$ and $k>j.$ Then $k>i$ and, therefore,
$a_{ik}>1.$ Any consistent matrix in $\mathcal{PC}_{n}$ is order transitive.

\bigskip

We next generalize the concept of order transitivity to block reciprocal
matrices. Let $S=\{S_{1},\ldots,S_{p}\}$ be a partition of $N$. We denote by
$\mathcal{PC}_{n}^{(S)}$ the set of matrices $A\in\mathcal{PC}_{n}$ such that,
for all $i,j\in\{1,\ldots,p\},$ $i\neq j$, either $A[S_{i},S_{j}]\geq J$ or
$A[S_{i},S_{j}]\leq J.$ Matrix $A\in\mathcal{PC}_{n}^{(S)}$ is said to be
$\{S_{1},\ldots,S_{p}\}$\emph{-order transitive} ($\{S_{1},\ldots,S_{p}\}$-OT,
for short) if, for all $i,j,k\in\{1,\ldots,p\},$ $A[S_{i},S_{j}]\geq J$ and
$A[S_{j},S_{k}]\geq J$ imply $A[S_{i},S_{k}]\geq J.$ If $i,j,k$ are not
distinct, the last condition is trivially satisfied. Note that $A[S_{i}%
,S_{j}]\geq J\ $implies that either $i\neq j$ or $A[S_{i},S_{j}]=J$.

\begin{example}
The matrix%
\[
A=\left[
\begin{tabular}
[c]{c|ccc|c}%
$1$ & $2$ & $1$ & $2$ & $2$\\\hline
$\frac{1}{2}$ & $1$ & $2$ & $\frac{1}{2}$ & $3$\\
$1$ & $\frac{1}{2}$ & $1$ & $2$ & $4$\\
$\frac{1}{2}$ & $2$ & $\frac{1}{2}$ & $1$ & $3$\\\hline
$\frac{1}{2}$ & $\frac{1}{3}$ & $\frac{1}{4}$ & $\frac{1}{3}$ & $1$%
\end{tabular}
\right]
\]
is $\{\{1\},\{2,3,4\},\{5\}\}$-OT. It is also is $\{\{1,2,3,4\},\{5\}\}$-OT
and is $\{\{1\},\{2,3,4,5\}\}$-OT. Note, however, that $A$ is not order
transitive, as $A[\{2,3,4\}]$ is not.
\end{example}

\bigskip

There exists a characterization of $\{S_{1},\ldots,S_{p}\}$-OT matrices that
extends the one mentioned above for order transitive matrices.

\begin{theorem}
\label{Prop2}Let $S=\{S_{1},\ldots,S_{p}\}$ be a partition of $N$ and let
$A\in\mathcal{PC}_{n}^{(S)}.$ If $A$ is $\{S_{1},\ldots,S_{p}\}$-OT, then
there is a permutation $i_{1}i_{2}\cdots i_{p}$ of $\{1,\ldots,p\}$ such that
$A[S_{i_{s}},S_{i_{t}}]\geq J$, for all $s,t\in\{1,\ldots,p\},$ $t>s.$ The
converse also holds if, for every $i,j\in\{1,\ldots,p\}$, $i\neq j$,
$A[S_{i},S_{j}]$ has an entry distinct from $1$.
\end{theorem}

\begin{proof}
To prove the first claim, suppose that $A$ is $\{S_{1},\ldots,S_{p}\}$-OT. Let
$T$ be the $p$-by-$p$ matrix with zero diagonal and such that, for $i\neq j,$
the $i,j$ entry is $1$ if $A[S_{i},S_{j}]\geq J$ and is $0$ otherwise. Since
either $A[S_{i},S_{j}]\geq J$ or $A[S_{i},S_{j}]\leq J$, it follows that $T$
is the adjacency matrix of a semicomplete digraph. By Theorem \ref{thRedeis},
there is an H-path $i_{1}i_{2}\cdots i_{p}$ in the digraph associated with
$T$. Thus, $A[S_{i_{s}},S_{i_{s+1}}]\geq J$ for all $s=1,\ldots,p-1.$ Now the
result follows since, by assumption, $A$ is $\{S_{1},\ldots,S_{p}\}$-OT.

To prove the converse claim, suppose that $A[S_{i_{j_{1}}},S_{i_{j_{2}}}]\geq
J$ and $A[S_{i_{j_{2}}},S_{i_{j_{3}}}]\geq J$, for some $j_{1},j_{2},j_{3}$
distinct. Since $A[S_{i_{s}},S_{i_{t}}]\geq J$ (and not $A[S_{i_{t}},S_{i_{s}%
}]\geq J$), for all $s,t=1,\ldots,p,$ $t>s$, it follows that $j_{2}>j_{1}$ and
$j_{3}>j_{2}.$ Then, $j_{3}>j_{1}$ implying $A[S_{i_{j_{1}}},S_{i_{j_{3}}%
}]\geq J$.
\end{proof}

\bigskip We observe that the fact that each block $A[S_{i},S_{j}]$, $i\neq j$,
has an entry different from $1$ is important in the proof of the converse
claim in Theorem \ref{Prop2}.

\bigskip If $S=\{S_{1},S_{2}\}$, then $A\in\mathcal{PC}_{n}^{(S)}$ is
$\{S_{1},S_{2}\}$-OT. When $p>2$, an $\{S_{1},\ldots,S_{p}\}$-OT matrix is
$\{S_{1}^{\prime},\ldots,S_{q}^{\prime}\}$-OT for a coarser partition
$\{S_{1}^{\prime},\ldots,S_{q}^{\prime}\}$. \bigskip

If $A\in\mathcal{PC}_{n}$ is order transitive, then $A$ is $\{S_{1}%
,\ldots,S_{p}\}$-OT for $S_{1}=\{1\},\ldots,S_{n}=\{n\}$. For such a matrix,
there is an efficient vector $w$ that has no rank reversals ($w$ may be
obtained from a certain spanning tree). This was shown in \cite{FJRR} when all
off-diagonal entries of $A$ are distinct from $1$, but the same argument can
be used otherwise. We next generalize this result to block order transitive matrices.

\begin{corollary}
\label{Cordertr}Let $S=\{S_{1},\ldots,S_{p}\}$ be a partition of $N$ and let
$A\in\mathcal{PC}_{n}^{(S)}$ be an $\{S_{1},\ldots,S_{p}\}$-OT matrix. Then,
there is a vector $w\in\mathcal{E}(A)$ such that $w[S_{i}]\in\mathcal{E}%
(A[S_{i}]),$ $i=1,\ldots,p$, and that has no rank reversals at positions
outside the principal blocks indexed by the sets $S_{i}$ in the partition.
\end{corollary}

\begin{proof}
Because of Theorem \ref{Prop2}, we may assume, without loss of generality,
that $A$ is $\{S_{1},\ldots,S_{p}\}$-OU for the path $12\cdots p$ with weights
$1$. Thus, $A[S_{j},S_{j+1}]\geq J$ for all $j\in\{1,\ldots,p-1\}$ and
$A[S_{j},S_{j+1}\cup\cdots\cup S_{p}]\geq J$ for all $j\in\{1,\ldots,p-1\}.$
Then the statement follows from Corollary \ref{cortrans}.
\end{proof}

\bigskip

An interesting situation is when any alternative in a set $S_{1}\subset N$ is
strictly ranked above any alternative not in $S_{1}$ by $A$ and any
alternative in $S_{3}\subset N\backslash S_{1}$ is strictly ranked below any
alternative not in $S_{3}$ by $A$. In this case all entries of $A$ that lie in
a row indexed by $S_{1}$ and in a column indexed by $N\backslash S_{1}$, or in
a column indexed by $S_{3}$ and in a row indexed by $N\backslash S_{3},$ are
$>1.$ We call such a reciprocal matrix an $S_{1}$\emph{-best, }$S_{3}%
$\emph{-worst matrix}. Such a matrix is permutationally similar to a block
matrix of the form%
\[
\left[
\begin{array}
[c]{ccc}%
A_{11} & A_{12} & A_{13}\\
A_{12}^{(-T)} & A_{22} & A_{23}\\
A_{13}^{(-T)} & A_{23}^{(-T)} & A_{33}%
\end{array}
\right]
\]
in which $A_{11}\in\mathcal{PC}_{\#S_{1}}$ and $A_{33}\in\mathcal{PC}%
_{\#S_{3}},$ $A_{12}>J,$ $A_{13}>J$ and $A_{23}>J$. Thus, $A$ is an
$\{S_{1},S_{2},S_{3}\}$-OT reciprocal matrix.

\begin{example}%
\[
\left[
\begin{tabular}
[c]{cc|cc|cc}%
$1$ & $2$ & $3$ & $2$ & $4$ & $2$\\
$\frac{1}{2}$ & $1$ & $5$ & $2$ & $3$ & $4$\\\hline
$\frac{1}{3}$ & $\frac{1}{5}$ & $1$ & $\frac{1}{4}$ & $6$ & $6$\\
$\frac{1}{2}$ & $\frac{1}{2}$ & $4$ & $1$ & $3$ & $3$\\\hline
$\frac{1}{4}$ & $\frac{1}{3}$ & $\frac{1}{6}$ & $\frac{1}{3}$ & $1$ & $2$\\
$\frac{1}{2}$ & $\frac{1}{4}$ & $\frac{1}{6}$ & $\frac{1}{3}$ & $\frac{1}{2}$
& $1$%
\end{tabular}
\right]
\]
is a $\{1,2\}$-best, $\{5,6\}$-worst matrix. Also, it is an
$\{\{1,2\},\{3,4\},\{5,6\}\}$-OT reciprocal matrix.
\end{example}

From Corollary \ref{Cordertr}, we get that the best/worst property is
preserved by some efficient vector.

\begin{corollary}
\label{T1}If $A\in\mathcal{PC}_{n}$ is an $S_{1}$-best, $S_{3}$-worst matrix,
then there is a $w\in\mathcal{E}(A)$ strictly ranking the alternatives in
$S_{1}$ above the remaining alternatives and the alternatives in $S_{3}$ below
the remaining alternatives.
\end{corollary}

The "strict" claim follows due to Theorem \ref{thstep}.

\section{Block uniformly ordered reciprocal matrices\label{SUO}}

Matrix $A\in\mathcal{PC}_{n}$ is said to be \emph{uniformly ordered} (UO) if
all vectors in $\mathcal{E}(A)$ have the same (weak) order. Such matrices were
introduced in \cite{FJ20}, where the case $n=4$ was studied in detail. Several
characterizations of UO reciprocal matrices $A$ using $L_{A}$ were given later
in \cite{FJ8,FJRR}. A useful one is that $A$ is UO if and only if there is a
matrix permutation similar to $L_{A}$ with all entries above the diagonal
$\geq1$ (equivalently, $L_{A}$ has $\frac{n^{2}-n}{2}$ off-diagonal entries
$\geq1$). If $A$ is UO, the same permutation similarity on $A$ gives a matrix
for which all efficient vectors are decreasing. In particular the permuted
matrix has decreasing columns. Here, we study a generalization of UO
reciprocal matrices to block reciprocal matrices.

\bigskip

Let $\{S_{1},\ldots,S_{p}\}$ be a partition of $N.$ Matrix $A\in
\mathcal{PC}_{n}$ is said to be $\{S_{1},\ldots,S_{p}\}$\emph{-uniformly
ordered} ($\{S_{1},\ldots,S_{p}\}$-UO, for short) if there is a permutation
$i_{1}i_{2}\cdots i_{p}$ of $\{1,\ldots,p\}$ such that, for all $w\in
\mathcal{E}(A)$ and all $j\in\{1,\ldots,p-1\}$, $W[S_{i_{j}},S_{i_{j+1}}]\geq
J$. Because $W$ is consistent, this implies that $W[S_{i_{s}},S_{i_{t}}]\geq
J$ for all $t>s$. Thus, for such an $A$, all efficient vectors give the unique
order $S_{i_{1}}\succsim\cdots\succsim S_{i_{p}}$ of the sets of alternatives.

\bigskip

In general an $\{S_{1},\ldots,S_{p}\}$-UO matrix is $\{S_{1}^{\prime}%
,\ldots,S_{q}^{\prime}\}$-UO for a coarser partition $\{S_{1}^{\prime}%
,\ldots,S_{q}^{\prime}\}$, just as with block order transitivity. If
$i_{1}i_{2}\cdots i_{p}$ as above is $12\cdots p,$ we may have $S_{1}^{\prime
}=\cup_{j=1}^{j_{1}}S_{j},$ $S_{2}^{\prime}=\cup_{j=j_{1}+1}^{j_{2}}S_{j},$
$\ldots,$ $S_{q}^{\prime}=\cup_{j=j_{q-1}+1}^{p}S_{j}$, with $1\leq
j_{1}<j_{2}<\cdots<j_{q-1}<p.$ A UO-matrix is $\{\{1\},\ldots,\{n\}\}$-UO.

\bigskip

A characterization similar to the one for UO matrices referred above is
possible for block UO matrices.

\begin{theorem}
\label{thUO}Let $A\in\mathcal{PC}_{n}^{0}$. Matrix $A$ is $\{S_{1}%
,\ldots,S_{p}\}$-UO if and only if there is a permutation $i_{1}i_{2}\cdots
i_{p}$ of $\{1,\ldots,p\}$ such that $L_{A}[S_{i_{j}},S_{i_{j+1}}]\geq J$ for
all $j=1,\ldots,p-1$.
\end{theorem}

\begin{proof}
The "only if" claim follows from Theorem \ref{rattain}. The "if" claim follows
from Corollary \ref{Thintgea}.
\end{proof}

\bigskip

\begin{remark}
\label{RR}Let $A\in\mathcal{PC}_{n}^{0}$. If $L_{A}[S_{i_{j}},S_{i_{j+1}}]\geq
J$ for all $j=1,\ldots,p-1$, then $L_{A}[S_{i_{s}},S_{i_{t}}]\geq J$ for all
$t>s.$ In fact, if the former condition occurs, by Corollary \ref{Thintgea},
for any $w\in\mathcal{E}(A)$ and any $j=1,\ldots,p-1$, we have $W[S_{i_{j}%
},S_{i_{j+1}}]\geq J.$ Since $W$ is consistent, this implies that $W[S_{i_{s}%
},S_{i_{t}}]\geq J$ for $t>s$. Then, the claim follows by Theorem
\ref{rattain}.
\end{remark}

We then have the following.

\begin{corollary}
\label{cper}If $A\in\mathcal{PC}_{n}^{0}$ is $\{S_{1},\ldots,S_{p}\}$-UO, then
there is a permutation $i_{1}i_{2}\cdots i_{p}$ of $\{1,\ldots,p\}$ such that
$A[S_{i_{s}},S_{i_{t}}]\geq J$, for all $t>s.$ Moreover, for all
$w\in\mathcal{E}(A)$, there are no rank reversals at positions outside the
principal blocks indexed by the sets $S_{i}$ in the partition.
\end{corollary}

\begin{proof}
The first claim follows from Theorem \ref{thUO} and Remark \ref{RR}, and the
fact that, by Corollary \ref{Thintgea}, $L_{A}\leq A$ (entrywise). Then the
second claim also follows from Corollary \ref{Thintgea}.
\end{proof}

\bigskip

Note that Corollary \ref{Thintgea} also allows us to conclude that, for any
column of $A$, the entries indexed by $S_{i_{j}}$ are higher than those
indexed by $S_{i_{j+1}},$ $j\in\{1,\ldots,p-1\}.$

\bigskip

If all blocks $A[S_{i},S_{j}]$, $i\neq j$, in an $\{S_{1},\ldots,S_{p}\}$-UO
matrix have an entry different from $1$, then by Theorem \ref{Prop2} and
Corollary \ref{cper}, $A$ is $\{S_{1},\ldots,S_{p}\}$-OT.

\begin{example}
Let
\[
A=\left[
\begin{tabular}
[c]{c|cc|c}%
$1$ & $2$ & $3$ & $20$\\\hline
$\frac{1}{2}$ & $1$ & $\frac{3}{4}$ & $5$\\
$\frac{1}{3}$ & $\frac{4}{3}$ & $1$ & $\frac{10}{3}$\\\hline
$\frac{1}{20}$ & $\frac{1}{5}$ & $\frac{3}{10}$ & $1$%
\end{tabular}
\right]  .
\]
We show that matrix $A$ is a $\{\{1\},\{2,3\},\{4\}\}$-UO matrix. We have
$\Gamma(A)=\{\alpha\},$ in which $\alpha=12341$. Moreover, $\alpha(A)=\frac
{1}{4}.$ Then
\[
L_{A}=P_{A,\alpha}=\left[
\begin{tabular}
[c]{c|cc|c}%
$1$ & $2$ & $\frac{3}{2}$ & $5$\\\hline
$\frac{1}{8}$ & $1$ & $\frac{3}{4}$ & $\frac{5}{2}$\\
$\frac{1}{6}$ & $\frac{1}{3}$ & $1$ & $\frac{10}{3}$\\\hline
$\frac{1}{20}$ & $\frac{1}{10}$ & $\frac{3}{40}$ & $1$%
\end{tabular}
\ \right]  .
\]

\end{example}

Thus, taking into account Corollary \ref{Thintgea}, all vectors in
$\mathcal{E}(A)$ give the strict order $\{1\}\succ\{2,3\}\succ\{4\}.$ Note
that $A$ is not UO (in the usual sense) as its columns do not have all the
same order.

\section{Conclusions \label{sconcl}}

Here we study efficient vectors for block reciprocal matrices. We have shown
how efficient vectors of individual principal blocks can be extended to
efficient vectors for the full matrix, generalizing an earlier extension
result in the scalar setting. We then investigated the preservation of ordinal
information between blocks. Conditions were obtained that guarantee the
existence of efficient vectors respecting the ordering implied by cross-block comparisons.

We introduced a notion of order transitivity for block reciprocal matrices. We
proved that, under order transitivity, there always exists an efficient vector
that avoids rank reversals associated with off-diagonal blocks.

The concept of uniformly ordered reciprocal matrices was also extended to the
block setting. We characterized the block uniformly ordered reciprocal
matrices through bounds for the consistent matrices generated by efficient
vectors. For these matrices, all efficient vectors lead to the same ordinal
ranking of the considered sets of alternatives, providing uniformity in the
decision outcome.

We obtained a graph-theoretic characterization for the existence of a positive
vector that is efficient for every completion of a given partial reciprocal
matrix. This result provides conditions guaranteeing the existence of an
efficient vector that preserves certain ordinal relations induced by
comparisons among blocks for all completions of the partial matrix.

Overall, the paper shows that several important notions from the theory of
reciprocal matrices can be naturally extended to block structures. These
extensions help clarify how local and global preference information interact
and provide tools for studying order, incomplete and large pairwise comparison problems.

\bigskip

\textbf{Declaration of conflict of interest}

\bigskip

There are no conflicts of interest.

\end{document}